\documentclass[12pt]{amsart}
\usepackage{amsxtra}
\usepackage{amssymb, amsmath}
\usepackage{amsthm}
\usepackage {hyperref}

\newtheorem{theorem}{Theorem}
\newtheorem{remark}{Remark}
\newtheorem{example}{Example}

\newtheorem{corollary}{Corollary}
\newtheorem{result}{Result}
\newtheorem{proposition}{Proposition}

\title[Jensen Inequality]{On a Difference of Jensen Inequality and its Applications to
Mean Divergence Measures}

\thanks{This paper is a part of author's chapter to appear in
\textit{'Advances in Imaging and Electron Physics'}, 2005,
Elsevier Publication}

\author{Inder Jeet Taneja}
\address{Inder Jeet Taneja\\
Departamento de Matem\'{a}tica\\ Universidade Federal
de Santa Catarina\\
88.040-900 Florian\'{o}polis, SC, Brazil}

\email{taneja@mtm.ufsc.br}

\urladdr{http://www.mtm.ufsc.br/$\sim $taneja}


\begin{document}

\begin{abstract}
In this paper we have considered a difference of Jensen's
inequality for convex functions and proved some of its properties.
In particular, we have obtained results for Csisz\'{a}r
\cite{csi1} $f-$divergence. A result is established that allow us
to compare two measures under certain conditions. By the
application of this result we have obtained a new inequality for
the well known means such as arithmetic, geometric and harmonic.
Some divergence measures based on these means are also defined.
\end{abstract}

\maketitle

\section{Jensen Difference}

Let
\[
\Gamma _n = \left\{ {P = (p_1 ,p_2 ,...,p_n )\left| {p_i >
0,\sum\limits_{i = 1}^n {p_i = 1} } \right.} \right\}, \,\, n
\geqslant 2,
\]

\noindent be the set of all complete finite discrete probability
distributions.

Let $f:I \subset \mathbb{R} \to \mathbb{R}$ be a differentiable
convex function on the interval $I$, $x_i \in \mathop I\limits^o $
($\mathop I\limits^o $ is the interior of $I)$. Let $\lambda =
(\lambda _1 ,\lambda _1 ,...,\lambda _n ) \in \Gamma _n $, then it
is well known that
\begin{equation}
\label{eq1} f\left( {\sum\limits_{i = 1}^n {\lambda_i x_i } }
\right) \leqslant \sum\limits_{i = 1}^n {\lambda_i f(x_i )} .
\end{equation}

The above inequality is famous as \textit{Jensen inequality}. If
$f$ is concave, the inequality sign changes.

Let us consider the following \textit{Jensen difference}:
\begin{equation}
\label{eq2} F_f (\lambda,X) = \sum\limits_{i = 1}^n {\lambda_i
f(x_i ) - f\left( {\sum\limits_{i = 1}^n {\lambda_i x_i } }
\right)} ,
\end{equation}

Here below we shall give two theorems giving properties of
\textit{Jensen difference}.

\begin{theorem} \label{the11} Let $f:I \subset \mathbb{R} \to
\mathbb{R}$ be a differentiable convex function on the interval
$I$, $x_i \in \mathop I\limits^o $ ($\mathop I\limits^o $ is the
interior of $I)$, $\lambda = (\lambda _1 ,\lambda _1 ,...,\lambda
_n ) \in \Gamma _n $. If $\eta _1 ,\,\,\eta _2 \in \mathop
I\limits^o $ and $\eta _1 \leqslant x_i \leqslant \eta _2 $,
$\forall i = 1,2,...,n$, then we have the inequalities:
\begin{equation}
\label{eq3} 0 \leqslant F_f (\lambda, X) \leqslant L_f (\lambda, X
) \leqslant Z_f (\eta _1 ,\eta _2 ),
\end{equation}

\noindent where
\begin{equation}
\label{eq4} L_f (\lambda, X ) = \sum\limits_{i = 1}^n {\lambda _i
x_i {f}'(x_i ) - \left( {\sum\limits_{i = 1}^n {\lambda _i x_i } }
\right)\left( {\sum\limits_{i = 1}^n {\lambda _i {f}'(x_i )} }
\right)}
\end{equation}

\noindent and
\begin{equation}
\label{eq5}
Z_f (\eta _1 ,\eta _2 ) = \frac{1}{4}(\eta _2 - \eta _1 )\left[ {{f}'(\eta
_2 ) - {f}'(\eta _1 )} \right].
\end{equation}
\end{theorem}

The above theorem is due to Dragomir \cite{dra3}. It has been
applied by many authors \cite{dra2},\cite{ddp}. The measure
$F(\lambda,X)$ has been extensively studied by Burbea and Rao
\cite{bur1, bur2}. As a consequence of above theorem we have the
following corollary.

\begin{corollary} \label{cor11} For all $a,b,\upsilon ,\omega \in
(0,\infty )$, the following inequalities hold:
\begin{align}
\label{eq6} 0 & \leqslant \frac{\upsilon f(a) + \omega
f(b)}{\upsilon + \omega } - f\left( {\frac{\upsilon a + \omega
b}{\upsilon + \omega }} \right)\\
& \leqslant \frac{\upsilon a{f}'(a) + \omega b{f}'(b)}{\upsilon +
\omega } - \left( {\frac{\upsilon a + \omega b}{\upsilon + \omega
}} \right)\left( {\frac{\upsilon {f}'(a) + \omega
{f}'(b)}{\upsilon + \omega }} \right)\notag\\
& \leqslant \frac{1}{4}(b - a)\left( {{f}'(b) - {f}'(a)}
\right).\notag
\end{align}
\end{corollary}

\begin{proof} It follows from Theorem \ref{the11}, by taking
$\lambda _1 = \frac{\upsilon }{\upsilon + \omega }$, $\lambda _2 =
\frac{\omega }{\upsilon + \omega }$, $\lambda _3 = ... = \lambda
_n = 0$, $x_1 = a$, $x_2 = b$, $x_2 = ... = x_n = 0$.
\end{proof}

Now we shall give some examples of Theorem \ref{the11}.

\begin{example} \label{exa11} For all $x \in (0,\infty )$, let us
consider a function
\begin{equation}
\label{eq7} f_s (x) = \begin{cases}
 {\frac{1 - x^s}{s},} & {s \ne 0,} \\
 { - \ln x,} & {s = 0.} \\
\end{cases}
\end{equation}

We can easily check that the function $f_s (x)$ is convex in
$(0,\infty )$ for all $s \leqslant 1$. Let there exist $\eta _1 $
and $\eta _2 $ such that $\eta _1 \leqslant x_i \leqslant \eta _2
$, $\forall i = 1,2,...,n$. Applying Theorem \ref{the11} for the
function $f_s (x)$, we have
\begin{equation}
\label{eq8} 0 \leqslant F_s (\lambda,X) \leqslant Z_s (\eta _1
,\eta _2 ), \,\, s \leqslant 1,
\end{equation}

\noindent where
\begin{equation}
\label{eq9} F_s (\lambda,X) = \begin{cases}
 {\frac{1}{s}\left[ \left(
{\sum\limits_{i = 1}^n {\lambda_i x_i } } \right)^s -
{\sum\limits_{i = 1}^n {\lambda_i x_i^s }} \right],} & {s \ne 0,}
\\\\
 {\ln \left( {\frac{A(\lambda,X)}{G(\lambda,X)}} \right),} & {s = 0.} \\
\end{cases}
\end{equation}
\begin{equation}
\label{eq10} A(\lambda,X) = \sum\limits_{i = 1}^n {\lambda_i x_i }
,
\end{equation}
\begin{equation}
\label{eq11} G(\lambda,X) = \prod\limits_{i = 1}^n {x_i^{\lambda_i
} }
\end{equation}

\noindent and
\begin{equation}
\label{eq12}
Z_s (\alpha ,\beta ) = \frac{1}{4}(\eta _2 - \eta _1 )\left( {\eta _1 ^{s -
1} - \eta _2 ^{s - 1}} \right).
\end{equation}
\end{example}

In particular we have
\begin{equation}
\label{eq13} \frac{A(\lambda,X)}{G(\lambda,X)} \leqslant \exp
\left[ {\frac{(\eta _2 - \eta _1 )^2}{4\eta _1 \eta _2 }}
\right],\,\,\eta _1 \leqslant x_i \leqslant \eta _2 ,\,\,\forall i
= 1,2,..n.
\end{equation}

\bigskip
The result (\ref{eq13}) is due to Dragomir \cite{dra3}. The
following proposition is a particular case of the inequalities
(\ref{eq6}) and gives bounds on Burbea and Rao's \cite{bur1, bur2}
\textit{Jensen difference divergence measure}.

\begin{proposition} \label{pro11} Let $f:(0,\infty ) \to \mathbb{R}$ be
a differentiable convex function. Then for all $P,Q \in \Gamma _n
$, we have
\begin{equation}
\label{eq14}
0 \leqslant \sum\limits_{i = 1}^n {\left[ {\frac{f(p_i ) + f(q_i )}{2} -
f\left( {\frac{p_i + q_i }{2}} \right)} \right]}
 \leqslant \frac{1}{4}\sum\limits_{i = 1}^n {\left( {p_i - q_i } \right)}
\left[ {{f}'(p_i ) - {f}'(q_i )} \right].
\end{equation}
\end{proposition}

\begin{proof} Take $\omega = \upsilon = \frac{1}{2}$ in
(\ref{eq6}), we get
\begin{equation} \label{eq15} 0 \leqslant \frac{f(a) + f(b)}{2} -
f\left( {\frac{a + b}{2}} \right) \,\,
 \leqslant \frac{1}{4}(b - a)\left[ {{f}'(b) - {f}'(a)} \right].
\end{equation}

Replace in (\ref{eq15}), $a$ by $p_i $ and $b$ by $q_i $, and sum
over all $i = 1,2,...,n$, we get the required result.
\end{proof}

\begin{example} \label{exa12} Let us consider a convex function
\begin{equation}
\label{eq16} \phi _s (x) = \begin{cases}
 {\left[ {s(s - 1)} \right]^{ - 1}\left[ {x^s - 1 - s(x - 1)} \right],} & {s
\ne 0,1}, \\
 {x - 1 - \ln x,} & {s = 0}, \\
 {1 - x + x\ln x,} & {s = 1}, \\
\end{cases}
\end{equation}

\noindent for all $x \in (0,\infty )$ and $s \in ( - \infty
,\infty )$. Then from (\ref{eq14}), we get
\begin{equation}
\label{eq17} 0 \leqslant \mathcal{W}_s (P\vert \vert Q) \leqslant
\frac{1}{4} \mathcal{V}_s (P\vert \vert Q),
\end{equation}

\noindent where
\begin{equation}
\label{eq18} \mathcal{W}_s (P\vert \vert Q) = \begin{cases}
 {I_s (P\vert \vert Q) = \left[ {s(s - 1)} \right]^{ - 1}\sum\limits_{i =
1}^n {\left[ {\frac{p_i^s + q_i^s }{2} - \left( {\frac{p_i + q_i }{2}}
\right)^s} \right],} } & {s \ne 0,1,} \\\\
 {I_0 (P\vert \vert Q) = \ln \left[ {\prod\limits_{i = 1}^n {\left(
{\frac{p_i + q_i }{2\sqrt {p_i q_i } }} \right)} } \right],} & {s = 0,}
\\\\
 {I(P\vert \vert Q) =  H\left( {\frac{P + Q}{2}}
\right) - \frac{H(P) + H(Q)}{2},} & {s = 1,} \\
\end{cases}
\end{equation}

\noindent and
\begin{equation}
\label{eq19} \mathcal{V}_s (P\vert \vert Q) = \begin{cases}
 {J_s (P\vert \vert Q) = \frac{1}{(s - 1)}\sum\limits_{i = 1}^n {\left(
{p_i - q_i } \right)\left( {p_i^{s - 1} - q_i^{s - 1} } \right),} } & {s \ne
0,1,} \\\\
 {J_0 (P\vert \vert Q) = \sum\limits_{i = 1}^n {\frac{\left( {p_i
- q_i } \right)^2}{p_i q_i }} ,} & {s = 0,} \\\\
 {J(P\vert \vert Q) = \sum\limits_{i = 1}^n {\left( {p_i - q_i }
\right)\ln \left( {\frac{p_i }{q_i }} \right)} ,} & {s = 1.} \\
\end{cases}
\end{equation}
\end{example}

The expression $H(P) = - \sum\limits_{i = 1}^n {p_i \ln p_i } $,
appearing in (\ref{eq18}) is the well known \textit{Shannon's
entropy}. The expression $J(P\vert \vert Q)$ appearing in
(\ref{eq19}) is \textit{Jeffreys-Kullback-Leibler's J-divergence}
(ref. Jeffreys \cite{jef} and Kullback and Leibler \cite{kul}).
The expression $J_s (P\vert \vert Q)$ is due to Burbea and Rao
\cite{bur1}. The measures (\ref{eq18}) and (\ref{eq19}) has been
studied by Burbea and Rao \cite{bur1} only for positive values of
the parameters. Some studies on these generalised measures can be
seen in Taneja \cite{tan2, tan4}. Here we have presented them for
all $s \in (-\infty, \infty)$. The function given in (\ref{eq16})
is due to Cressie and Read \cite{crr}.

\begin{proposition} \label{pro12} Let $f:\mathbb{R}_ + \to \mathbb{R}$
be differentiable convex and normalized i.e., $f(1) = 0$. If $P,Q
\in \Gamma _n $, are such that $0 < r \leqslant \frac{p_i }{q_i }
\leqslant R < \infty $, $\forall i \in \{1,2,...,n\}$, for some
$r$ and $R$ with $0 < r \leqslant 1 \leqslant R < \infty $, then
we have
\begin{equation}
\label{eq20}
0 \leqslant C_f (P\vert \vert Q) \leqslant E_{C_f } (P\vert \vert Q)
\leqslant A_{C_f } (r,R),
\end{equation}

\noindent and
\begin{equation}
\label{eq21}
0 \leqslant C_f (P\vert \vert Q) \leqslant B_{C_f } (r,R) \leqslant A_{C_f }
(r,R),
\end{equation}

\noindent where
\begin{equation}
\label{eq22}
C_f (P\vert \vert Q) = \sum\limits_{i = 1}^n {q_i } f(\frac{p_i }{q_i }),
\end{equation}
\begin{equation}
\label{eq23}
E_{C_f } (P\vert \vert Q) = \sum\limits_{i = 1}^n {(p_i - q_i )}
{f}'(\frac{p_i }{q_i }),
\end{equation}
\begin{equation}
\label{eq24} A_{C_f } (r,R) = \frac{1}{4}(R - r)\left( {{f}'(R) -
{f}'(r)} \right)
\end{equation}

\noindent and
\begin{equation}
\label{eq25}
B_{C_f } (r,R) = \frac{(R - 1)f(r) + (1 - r)f(R)}{R - r}.
\end{equation}
\end{proposition}

The inequalities (\ref{eq20}) follow in view of (\ref{eq3}). The
inequalities (\ref{eq21}) follow in view of (\ref{eq6}). For
details refer to Taneja \cite{tan8}. The above proposition is an
improvement over the work of Dragomir \cite{dra4, dra5}. The
measure (\ref{eq22}) is known as \textit{Csisz\'{a}r's \cite{csi1}
f-divergence}.

\begin{example} \label{exa13} Under the conditions of Proposition \ref{pro12},
the inequalities (\ref{eq20}) and (\ref{eq21}) for the function
(\ref{eq16}) are given by
\begin{equation}
\label{eq26}
0 \leqslant \Phi _s (P\vert \vert Q) \leqslant E_{\Phi _s } (P\vert \vert Q)
\leqslant A_{\Phi _s } (r,R)
\end{equation}

\noindent and
\begin{equation}
\label{eq27}
0 \leqslant \Phi _s (P\vert \vert Q) \leqslant B_{\Phi _s } (r,R) \leqslant
A_{\Phi _s } (r,R),
\end{equation}

\noindent where
\begin{equation}
\label{eq28} \Phi _s (P\vert \vert Q) = \begin{cases}
 {K_s (P\vert \vert Q) = \left[ {s(s - 1)} \right]^{ - 1}\left[
{\sum\limits_{i = 1}^n {p_i^s q_i^{1 - s} } - 1} \right],} & {s \ne 0,1,}
\\\\
 {K(Q\vert \vert P) = \sum\limits_{i = 1}^n {q_i \ln \left( {\frac{q_i }{p_i
}} \right)} ,} & {s = 0,} \\\\
 {K(P\vert \vert Q) = \sum\limits_{i = 1}^n {p_i \ln \left( {\frac{p_i }{q_i
}} \right)} ,} & {s = 1,} \\
\end{cases}
\end{equation}
\begin{equation}
\label{eq29} E_{\Phi _s } (P\vert \vert Q) = \begin{cases}
 {(s - 1)^{ - 1}\sum\limits_{i = 1}^n {(p_i - q_i )\left( {\frac{p_i }{q_i
}} \right)^{s - 1},} } & {s \ne 1,} \\\\
 {\sum\limits_{i = 1}^n {(p_i - q_i )\ln \left( {\frac{p_i }{q_i }}
\right),} } & {s = 1,} \\
\end{cases}
\end{equation}
\begin{equation}
\label{eq30} A_{\Phi _s } (r,R) = \frac{1}{4}\begin{cases}
 {\frac{(R - r)\left( {R^{s - 1} - r^{s - 1}} \right)}{4(s - 1)},} & {s \ne
1,} \\\\
 {\frac{1}{4}(R - r)\ln \left( {\frac{R}{r}} \right),} & {s = 1,} \\
\end{cases}
\end{equation}

\noindent and
\begin{equation}
\label{eq31} B_{\Phi _s } (r,R) = \begin{cases}
 {\frac{(R - 1)(r^s - 1) + (1 - r)(R^s - 1)}{(R - r)s(s - 1)},} & {s \ne
0,1,} \\\\
 {\frac{(R - 1)\ln \frac{1}{r} + (1 - r)\ln \frac{1}{R}}{(R - r)},} & {s =
0,} \\\\
 {\frac{(R - 1)r\ln r + (1 - r)R\ln R}{(R - r)},} & {s = 1.} \\
\end{cases}
\end{equation}
\end{example}

The measure $K(P\vert \vert Q)$ appearing in (\ref{eq28}) is the
well known \textit{Kullback-Leibler's \cite{kul} relative
information}. The measure $\Phi _s (P\vert \vert Q)$ given in
(\ref{eq28}) has been extensively studied in \cite{tan5},
\cite{tak}.

\begin{theorem} \label{the12} Let $f_1 ,f_2 :[a,b] \subset \mathbb{R}_
+ \to \mathbb{R}$ be twice differentiable functions on $(a,b)$ and
there are $\alpha $ and $\beta $ such that
\begin{equation}
\label{eq32} \alpha \leqslant \frac{{f}''_1 (x)}{{f}''_2 (x)}
\leqslant \beta , \,\, \forall x \in (a,b), \,\, {f}''_2 (x) > 0
\end{equation}

If $x_i \in [a,b]$ and $\lambda = (\lambda _1 ,\lambda _2
,...,\lambda _n ) \in \Gamma _n $, then
\begin{equation}
\label{eq33} \alpha \,F_{f_2 } (\lambda ,X) \leqslant F_{f_1 }
(\lambda ,X) \leqslant \beta \, F_{f_2 } (\lambda ,X).
\end{equation}
\begin{align}
\label{eq34} \alpha \left[ {L_{f_2 } (\lambda ,X) - F_{f_2 }
(\lambda ,X) } \right] & \leqslant L_{f_1 }
(\lambda ,X) - F_{f_1 } (\lambda ,X) \\
& \leqslant \beta \left[ L_{f_2 } (\lambda ,X) - {F_{f_2 }
(\lambda ,X) } \right]\notag
\end{align}

\noindent and
\begin{align}
\label{eq35} \alpha \left[ Z_{f_2 } (\eta _1 ,\eta _2 ) - {F_{f_2
} (\lambda ,X) } \right] & \leqslant Z_{f_1} (\eta _1 ,\eta _2 ) -
F_{f_1 } (\lambda ,X)
\\
& \leqslant \beta \left[ Z_{f_2 } (\eta _1 ,\eta _2 ) - {F_{f_2 }
(\lambda ,X) } \right].\notag
\end{align}
\end{theorem}

\begin{proof} Consider the mapping $g:[a,b] \to \mathbb{R}$,
defined by
\begin{equation}
\label{eq36} g(x) = f_1 (x) - \alpha f_2 (x), \quad \forall x \in
[a,b],
\end{equation}

\noindent where the functions $f_1 $ and $f_2 $ satisfy the
condition (\ref{eq32}). Then the function $g$ is twice
differentiable on $(a,b)$. This gives
\[
{g}'(x) = f_1 ^\prime (x) - \alpha f_2 ^\prime (x)
\]

\noindent and
\[
{g}''(x) = f_1 ^{\prime \prime }(x) - \alpha f_2 ^{\prime \prime
}(x) = f_2 ^{\prime \prime }(x)\left( {\frac{f_1 ^{\prime \prime
}(x)}{f_2 ^{\prime \prime }(x)} - \alpha } \right) \geqslant 0,
\,\,  \forall x \in (a,b).
\]

The above expression shows that $g$ is convex on $[a,b]$. Applying
Jensen inequality for the convex function $g$ one gets
\[
g\left( {\sum\limits_{i = 1}^n {\lambda _i x_i } } \right)
\leqslant \sum\limits_{i = 1}^n {\lambda _i g(x_i )},
\]

\noindent i.e.,
\[
f_1 \left( {\sum\limits_{i = 1}^n {\lambda _i x_i } } \right) -
\alpha f_2 \left( {\sum\limits_{i = 1}^n {\lambda _i x_i } }
\right) \leqslant \sum\limits_{i = 1}^n {\lambda _i \left[ {f_1
(x_i ) - \alpha f_2 (x)} \right]},
\]

\noindent i.e.,
\begin{equation}
\label{eq37} \alpha \left[ {\sum\limits_{i = 1}^n {\lambda _i f_2
(x_i ) - f_2 \left( {\sum\limits_{i = 1}^n {\lambda _i x_i } }
\right)} } \right] \leqslant \sum\limits_{i = 1}^n {\lambda _i f_1
(x_i ) - f_1 \left( {\sum\limits_{i = 1}^n {\lambda _i x_i } }
\right)} .
\end{equation}

The expression (\ref{eq37}) gives the $l.h.s.$ of the inequalities
(\ref{eq33}).

Again consider the mapping $k:[a,b] \to \mathbb{R}$ given by
\begin{equation}
\label{eq38} k(x) = \beta f_2 (x) - f_1 (x),
\end{equation}

\noindent and proceeding on similar lines as before, we get the
proof of the $r.h.s.$ of the inequalities (\ref{eq33}).

Now we shall prove the inequalities (\ref{eq34}). Applying the
inequalities (\ref{eq3}) for the convex function $g$ given by
(\ref{eq36}), we get
\[
F_g (\lambda ,X) \leqslant L_g (\lambda ,X) \leqslant Z_g (\eta _1
,\eta _2 ).
\]

\noindent i.e.,
\begin{align}
\label{eq39} F_{f_1 } (\lambda ,X) - \alpha F_{f_2 } (\lambda ,X)
& \leqslant L_{f_1 } (\lambda ,X) - \alpha L_{f_2 } (\lambda ,X)\\
& \leqslant Z_{f_1 } (\lambda ,X) - \alpha F_{f_2 } (\eta _1 ,\eta
_2 ).\notag
\end{align}

Simplifying the first inequality of (\ref{eq39}) we get
\begin{equation}
\label{eq40} \alpha \left[ {L_{f_2 } (\lambda ,X) - F_{f_2 }
(\lambda ,X)} \right] \leqslant L_{f_1 } (\lambda ,X) - F_{f_1 }
(\lambda ,X).
\end{equation}

Again simplifying the last inequality of (\ref{eq39}) we get
\begin{equation}
\label{eq41} \alpha \left[ {Z_{f_2 } (\eta _1 ,\eta _2 ) - F_{f_2
} (\lambda ,X)} \right] \leqslant Z_{f_1 } (\eta _1 ,\eta _2 ) -
F_{f_1 } (\lambda ,X).
\end{equation}

The expressions (\ref{eq40}) and (\ref{eq41}) complete the first
part of the inequalities (\ref{eq34}) and (\ref{eq35})
respectively. The last part of the inequalities (\ref{eq34}) and
(\ref{eq35}) follows by considering the function $k(x)$ given by
(\ref{eq38}) over the inequalities (\ref{eq3}).
\end{proof}

Particular cases of above theorem can be seen in \cite{anr},
\cite{dra1}, \cite{dra2}, \cite{drn}. Applications of the above
theorem for the \textit{Csisz\'{a}r's f-divergence} are given in
the following proposition.

\begin{proposition} \label{pro13} Let $f_1 ,f_2 :I \subset \mathbb{R}_
+ \to \mathbb{R}$ be two normalized convex mappings, i.e., $f_1
(1) = f_2 (1) = 0$ and suppose the assumptions:

(i) $f_1 $ and $f_2 $ are twice differentiable on $(r,R)$, where $0 < r
\leqslant 1 \leqslant R < \infty $;

(ii) there exists the real constants $\alpha ,\beta $ such that
$\alpha < \beta $ and
\begin{equation}
\label{eq42} \alpha \leqslant \frac{f_1 ^{\prime \prime }(x)}{f_2
^{\prime \prime }(x)} \leqslant \beta , \,\, f_2 ^{\prime \prime
}(x) > 0, \,\, \forall x \in (r,R).
\end{equation}

If $P,Q \in \Gamma _n $ are discrete probability distributions
satisfying the assumption
\[
0 < r \leqslant \frac{p_i }{q_i } \leqslant R < \infty ,
\]

\noindent then we have the inequalities:
\begin{equation}
\label{eq43} \alpha \,\,C_{f_2 } (P\vert \vert Q) \leqslant C_{f_1
} (P\vert \vert Q) \leqslant \beta \,\,C_{f_2 } (P\vert \vert Q),
\end{equation}
\begin{align}
\label{eq44} \alpha \left[ {E_{f_2 } (P\vert \vert Q) - C_{f_2 }
(P\vert \vert Q)} \right] & \leqslant E_{f_1 } (P\vert \vert Q) -
C_{f_1 } (P\vert \vert Q)\\
& \leqslant \beta \left[ {E_{f_2 } (P\vert \vert Q) - C_{f_2 }
(P\vert \vert Q)} \right],\notag
\end{align}
\begin{align}
\label{eq45} \alpha \left[ {A_{f_2 } (r,R) - C_{f_2 } (P\vert
\vert Q)} \right] & \leqslant A_{f_1 } (r,R) - C_{f_1 } (P\vert
\vert Q)\\
& \leqslant \beta \left[ {A_{f_2 } (r,R) - C_{f_2 } (P\vert \vert
Q)} \right]\notag
\end{align}

\noindent and
\begin{align}
\label{eq46} \alpha \left[ {B_{f_2 } (r,R) - C_{f_2 } (P\vert
\vert Q)} \right] & \leqslant B_{f_1 } (r,R) - C_{f_1 } (P\vert
\vert Q)\\
& \leqslant \beta \left[ {B_{f_2 } (r,R) - C_{f_2 } (P\vert \vert
Q)} \right].\notag
\end{align}
\end{proposition}

\begin{proof} It is an immediate consequence of the Theorem
\ref{the12}.
\end{proof}

\section{Applications to Mean Divergence Measures}

Let us consider the following mean of order $t$:
\begin{equation}
\label{eq47} D_t (a,b) = \begin{cases}
 {\left( {\frac{a^t + b^t}{2}} \right)^{1 / t},} & {t \ne 0}, \\
 {\sqrt {ab} ,} & {t = 0}, \\
 {\max \{a,b\},} & {t = \infty }, \\
 {\min \{a,b\},} & {t = - \infty } ,\\
\end{cases}
\end{equation}

\noindent for all $a,b >0$ and $t \in \mathbb{R}$. In particular,
we have

\begin{align}
D_{ - 1} (a,b) & = H(a,b) = \frac{2}{\frac{1}{a} + \frac{1}{b}}
= \frac{2ab}{a + b} = A(a^{ - 1},b^{ - 1})^{ - 1},\notag\\
D_0 (a,b) & = G(a,b) = \sqrt {ab} = \sqrt {A(a,b)H(a,b)}
,\notag\\
D_{1 / 2} (a,b) & = N_1 (a,b) = \left( {\frac{\sqrt a + \sqrt b
}{2}} \right)^2 = A\left( {\sqrt a ,\sqrt b } \right)^2\notag\\
\intertext{and} D_1 (a,b) & = A(a,b) = \frac{a + b}{2},\notag
\end{align}

\noindent where $H(a,b)$, $G(a,b)$ and $A(a,b)$ are the well known
\textit{harmonic}, \textit{geometric} and \textit{arithmetic
means} respectively. It is well know \cite{beb} that the
\textit{mean of order }$t$ given in (\ref{eq47}) is monotonically
increasing in $t$, then we can write
\[
D_{ - 1} (a,b) \leqslant D_0 (a,b) \leqslant D_{1 / 2} (a,b)
\leqslant D_1 (a,b),
\]

\noindent or equivalently,
\begin{equation}
\label{eq48}
H(a,b) \leqslant G(a,b) \leqslant N_1 (a,b) \leqslant A(a,b).
\end{equation}

We can easily check that the function $f(x) = - x^{1 / 2}$ is
convex in $(0,\infty )$. This allows us to conclude the following
inequality:
\begin{equation}
\label{eq49}
\frac{\sqrt a + \sqrt b }{2} \leqslant \sqrt {\frac{a + b}{2}} .
\end{equation}

From (\ref{eq49}), we can easily derive that
\begin{equation}
\label{eq50}
\left( {\frac{\sqrt a + \sqrt b }{2}} \right)^2 \leqslant \left(
{\frac{\sqrt a + \sqrt b }{2}} \right)\left( {\sqrt {\frac{a + b}{2}} }
\right) \leqslant \frac{a + b}{2}.
\end{equation}

Finally, the expressions (\ref{eq48}) and (\ref{eq50}) lead us to
following inequalities:
\begin{equation}
\label{eq51}
H(a,b) \leqslant G(a,b) \leqslant N_1 (a,b) \leqslant N_2 (a,b) \leqslant
A(a,b),
\end{equation}

\noindent where
\[
N_2 (a,b) = \left( {\frac{\sqrt a + \sqrt b }{2}} \right)\left( {\sqrt
{\frac{a + b}{2}} } \right).
\]

Let $P,Q \in \Gamma _n $. In (\ref{eq51}), replace $a$ by $p_i $
and $b$ by $q_i $ sum over all $i = 1,2,...n$ we get
\begin{equation}
\label{eq52} H(P\vert \vert Q) \leqslant G(P\vert \vert Q)
\leqslant N_1 (P\vert \vert Q) \leqslant N_2 (P\vert \vert Q)
\leqslant 1.
\end{equation}

Based on inequalities (\ref{eq52}), we shall build some
\textit{mean divergence measures}. Let us consider the following
differences:
\begin{align}
\label{eq53} M_{AG} (P\vert \vert Q) & = 1 - G(P\vert \vert Q),\\
\label{eq54} M_{AH} (P\vert \vert Q) & = 1 - H(P\vert \vert Q),\\
\label{eq55} M_{AN_2 } (P\vert \vert Q) & = 1 - N_2 (P\vert \vert
Q),\\
\label{eq56} M_{N_2 G} (P\vert \vert Q) & = N_2 (P\vert \vert Q) -
G(P\vert \vert Q),\\
\intertext{and}
\label{eq57} M_{N_2 N_1 } (P\vert \vert Q) & = N_2
(P\vert \vert Q) - N_1 (P\vert \vert Q).
\end{align}

We can easily verify that
\begin{align}
\label{eq58} M_{AG} (P\vert \vert Q) & = 1 - G(P\vert \vert Q)\\
 & = 2\left[ {N_1 (P\vert \vert Q) - G(P\vert \vert Q)} \right]: = 2M_{N_1 G}
(P\vert \vert Q)\notag\\
& = 2\left[ {1 - N_1 (P\vert \vert Q)} \right]: = 2M_{AN_1 }
(P\vert \vert Q).\notag\\
\end{align}

We can also write
\begin{equation}
\label{eq59}
M_{AG} (P\vert \vert Q) = 1 - G(P\vert \vert Q): = h(P\vert \vert Q)
\end{equation}

\noindent and
\begin{equation}
\label{eq60}
M_{AH} (P\vert \vert Q) = 1 - H(P\vert \vert Q): = \frac{1}{2}\Delta (P\vert
\vert Q),
\end{equation}

\noindent where $h(P\vert \vert Q)$ is the well known
\textit{Hellinger's \cite{hel} discrimination} and $\Delta (P\vert
\vert Q)$ is known by \textit{triangular discrimination}. These
two measures are well known in the literature of statistics. The
measure $M_{AN_2 } (P\vert \vert Q)$ is new and has been recently
studied by Taneja \cite{tan8}.\\

Now we shall prove the convexity of these measures. This is based
on the well known result due to Csisz\'{a}r \cite{csi1, csi2}.

\begin{result} \label{res21} If the function $f:\mathbb{R}_ + \to
\mathbb{R}$ is convex and normalized, i.e., $f(1) = 0$, then the
\textit{f-divergence}, $C_f (P\vert \vert Q)$ is
\textit{nonnegative} and \textit{convex} in the pair of
probability distribution $(P,Q) \in \Gamma _n \times \Gamma _n $.
\end{result}

\begin{example} \label{exa21} Let us consider
\begin{equation}
\label{eq61} f_{AH} (x) = \frac{(x - 1)^2}{2(x + 1)}, \,\, x \in
(0,\infty ),
\end{equation}

\noindent in (\ref{eq15}), then $C_f (P\vert \vert Q) = M_{AH}
(P\vert \vert Q),$ where $M_{AH} (P\vert \vert Q)$ is as given by
(\ref{eq54}).

Moreover,
\[
{f}'_{AH} (x) = \frac{(x - 1)(x + 3)}{2(x + 1)^2}
\]

\noindent and
\begin{equation}
\label{eq62} {f}''_{AH} (x) = \frac{4}{(x + 1)^3} > 0,\,\,x \in
(0,\infty ).
\end{equation}
\end{example}

\begin{example} \label{exa22} Let us consider
\begin{equation}
\label{eq63} f_{AG} (x) = \frac{1}{2}(\sqrt x - 1)^2, \,\, x \in
(0,\infty ),
\end{equation}

\noindent in (\ref{eq15}), then $C_f (P\vert \vert Q) = M_{AG}
(P\vert \vert Q),$ where $M_{AG} (P\vert \vert Q)$ is as given by
(\ref{eq53}).

Moreover,
\[
{f}'_{AG} (x) = \frac{\sqrt x - 1}{2\sqrt x }
\]

\noindent and
\begin{equation}
\label{eq64} {f}''_{AG} (x) = \frac{1}{4x\sqrt x } > 0,\,\,x \in
(0,\infty ).
\end{equation}
\end{example}

\begin{example} \label{exa23} Let us consider
\begin{equation}
\label{eq65} f_{N_2 N_1 } (x) = \frac{(x + 1)\sqrt {2(x + 1)} - 1
- x - 2\sqrt x }{4}, \,\, x \in (0,\infty )
\end{equation}

\noindent in (\ref{eq15}), then we have $C_f (P\vert \vert Q) =
M_{N_2 N_1 } (P\vert \vert Q),$ where $M_{N_2 N_1 } (P\vert \vert
Q)$ is as given by (\ref{eq57}).

Moreover,
\[
{f}'_{N_2 N_1 } (x) = \frac{2x + 1 + \sqrt x - \left( {\sqrt x + 1}
\right)\sqrt {2(x + 1)} }{6\sqrt x (x + 1)^2}
\]

\noindent and
\begin{align}
\label{eq66} {f}''_{N_2 N_1 } (x) & = \frac{ - 2x - 2x^{5 / 2} +
x(2x + 2)^{3 / 2}}{8x^{5 / 2}(2x + 2)^{3 / 2}}\\
& = \frac{x\left[ {(2x + 2)^{3 / 2} - 2(x^{3 / 2} + 1)}
\right]}{8x^{5 / 2}(2x + 2)^{3 / 2}}.\notag
\end{align}

Since $(x + 1)^{3 / 2} \geqslant x^{3 / 2} + 1$, $\forall x \in
(0,\infty )$ and $2^{3 / 2} \geqslant 2$, then obviously,
${f}''_{N_2 N_1 } (x) \geqslant 0$, $\forall x \in (0,\infty )$.
\end{example}

\begin{example} \label{exa24} Let us consider
\begin{equation}
\label{eq67} f_{N_2 G} (x) = \frac{\left( {\sqrt x + 1}
\right)\sqrt {2(x + 1)} - 4x}{4}, \,\, x \in (0,\infty ),
\end{equation}

\noindent in (\ref{eq15}), then $C_f (P\vert \vert Q) = M_{N_2 G}
(P\vert \vert Q),$ where $M_{N_2 G} (P\vert \vert Q)$ is as given
by (\ref{eq56}).

Moreover,
\[
{f}'_{_{N_2 G} } (x) = \frac{2x + 1 + \sqrt x - 2\sqrt {2(x + 1)} }{4\sqrt
{2x(x + 1)} }
\]

\noindent and
\begin{equation}
\label{eq68}
{f}''_{N_2 G} (x) = \frac{(2x + 2)^{3 / 2} - x^{3 / 2} - 1}{4x^{3 / 2}(2x +
2)^{3 / 2}}.
\end{equation}

Since $(x + 1)^{3 / 2} \geqslant x^{3 / 2} + 1$, $\forall x \in
(0,\infty )$ and $2^{3 / 2} \geqslant 1$, then obviously,
${f}''_{N_2 G} (x) \geqslant 0$, $\forall x \in (0,\infty )$.
\end{example}

\begin{example} \label{exa25} Let us consider
\begin{equation}
\label{eq69} f_{AN_2 } (x) = \frac{2(x + 1) - \left( {\sqrt x + 1}
\right)\sqrt {2(x + 1)} }{4}, \,\, x \in (0,\infty ),
\end{equation}

\noindent in (3.1), then $C_f (P\vert \vert Q) = M_{AN_2 } (P\vert
\vert Q),$ where $M_{AN_2 } (P\vert \vert Q)$ is as given by
(\ref{eq55}).

Moreover,
\[
{f}'_{AN_2 } (x) = - \frac{2x + 1 + \sqrt x - 2\sqrt {2x(x + 1)} }{4\sqrt
{2(x + 1)} },
\]

\noindent and
\begin{equation}
\label{eq70} {f}''_{AN_2 } (x) = \frac{1 + x^{3 / 2}}{8x^{3 / 2}(x
+ 1)\sqrt {2x + 2} } > 0,\,\,x \in (0,\infty ).
\end{equation}
\end{example}

In the above examples \ref{exa21}-\ref{exa25} the generating
function $f_{( \cdot )} (1) = 0$ and the second derivative is
positive for all $x \in (0,\infty )$. This proves the
\textit{nonegativity} and \textit{convexity} of the measures
(\ref{eq53})-(\ref{eq57}) in the pair of probability distributions
$(P,Q) \in \Gamma _n \times \Gamma _n $.

The inequality (\ref{eq52}) also admits more nonnegative
differences, but here we have considered only the convex ones.

Based on the Proposition \ref{pro12}, we can obtain bounds on the
\textit{mean divergence measures}, but we omit these details here.
Now we shall apply the inequalities (\ref{eq34}) given in
Proposition 1.3 to obtain inequalities among the measures
(\ref{eq53})-(\ref{eq57}).

\begin{theorem} \label{the21} The following inequalities among the six
\textit{mean divergences} hold:
\begin{equation}
\label{eq71}
\frac{1}{8}M_{AH} (P\vert \vert Q) \leqslant M_{N_2 N_1 } (P\vert \vert Q)
\leqslant \frac{1}{3}M_{N_2 G} (P\vert \vert Q)
\end{equation}
\[
 \leqslant \frac{1}{4}M_{AG} (P\vert \vert Q) \leqslant M_{AN_2 } (P\vert
\vert Q).
\]
\end{theorem}

The proof of the above theorem is based on the following
propositions, where we have proved each part separately.

\begin{proposition} \label{pro21} The following inequality hold:
\begin{equation}
\label{eq72} \frac{1}{8}M_{AH} (P\vert \vert Q) \leqslant M_{N_2
N_1 } (P\vert \vert Q).
\end{equation}
\end{proposition}

\begin{proof} Let us consider
\begin{align}
\label{eq73} g_{AH\_N_2 N_1 } (x) & = \frac{{f}''_{AH}
(x)}{{f}''_{N_2 N_1 } (x)}\\
& = \frac{32x^{5 / 2}(2x + 2)^{3 / 2}}{(x + 1)^3\left[ { - 2x -
2x^{5 / 2} + x(2x + 2)^{3 / 2}} \right]}, \,\, x \in (0,\infty
),\notag
\end{align}

\noindent where ${f}''_{AH} (x)$ and ${f}''_{N_2 N_1 } (x)$ are as
given by (\ref{eq62}) and (\ref{eq66}) respectively.

From (\ref{eq73}), we have
\begin{align}
{g}'_{AH\_N_2 N_1 } (x) & = - \frac{48\sqrt {2x(x + 1)} }{(x +
1)^4\left[ { - 2x - 2x^{5 / 2} + x(2x + 2)^{3 / 2}}
\right]^2}\times\notag\\
& \qquad \times \left[ {4x^2(1 - x^{5 / 2}) + x^2(x - 1)(2x +
2)^{5
/ 2}} \right]\notag \\
& = \frac{48x^2(x + 1)\left( {1 - \sqrt x } \right)\sqrt {2x(x +
1)} }{(x + 1)^4\left[ { - 2x - 2x^{5 / 2} + x(2x + 2)^{3 / 2}}
\right]^2}\times\notag\\
& \qquad \times \left[ {\sqrt 2 \left( {\sqrt x + 1} \right)\left(
{x + 1} \right)^{3 / 2} - \left( {x^2 + x^{3 / 2} + x + \sqrt x +
1} \right)} \right].\notag
\end{align}

Since $\sqrt {2(x + 1)} \geqslant \sqrt x + 1$, $\forall x \in
(0,\infty )$, then
\begin{align}
\sqrt 2 (x + 1)^{3 / 2}\left( {\sqrt x + 1} \right) & \geqslant
\left( {\sqrt x + 1} \right)^2(x + 1)\notag\\
& \geqslant x^2 + x^{3 / 2} + x + \sqrt x + 1.\notag
\end{align}

Thus we conclude that
\begin{equation}
\label{eq74} {g}'_{AH\_N_2 N_1 } (x)\begin{cases}
 { < 0,} & {x > 1}, \\
 { > 0,} & {x < 1}. \\
\end{cases}
\end{equation}

In view of (\ref{eq74}), we conclude that the function $g_{AH\_N_2
N_1 } (x)$ is increasing in $x \in (0,1)$ and decreasing in $x \in
(1,\infty )$, and hence
\begin{equation}
\label{eq75}
M = \mathop {\sup }\limits_{x \in (0,\infty )} g_{AH\_N_2 N_1 } (x) =
g_{AH\_N_2 N_1 } (1) = 8.
\end{equation}

Applying the inequalities (\ref{eq34}) for the measures $M_{AH}
(P\vert \vert Q)$ and $M_{N_2 N_1 } (P\vert \vert Q)$ along with
(\ref{eq75}) we get the required result.
\end{proof}

\begin{proposition} \label{pro22} The following inequality hold:
\begin{equation}
\label{eq76}
M_{N_2 N_1 } (P\vert \vert Q) \leqslant \frac{1}{3}M_{N_2 G} (P\vert \vert
Q).
\end{equation}
\end{proposition}

\begin{proof} Let us consider
\begin{align}
\label{eq77} g_{N_2 N_1 \_N_2 G} (x) & = \frac{{f}''_{N_2 N_1 }
(x)}{{f}''_{N_2 G} (x)} \\
& = \frac{ - 2x - 2x^{5 / 2} + x(2x + 2)^{3 / 2}}{2x\left[ {1 +
x^{3 / 2} - (2x + 2)^{3 / 2}} \right]}, \,\, x \in (0,\infty
),\notag
\end{align}

\noindent where ${f}''_{N_2 N_1 } (x)$ and ${f}''_{N_2 G} (x)$ are
as given by (\ref{eq66}) and (\ref{eq68}) respectively.

From (\ref{eq77}), we have
\begin{equation}
\label{eq78}
{g}'_{N_2 N_1 \_N_2 G_1 } (x) = \frac{3x^2\sqrt {2x + 2} \left( {1 - \sqrt x
} \right)}{2x^2\left[ { - 1 - x^{3 / 2} + (2x + 2)^{3 / 2}} \right]^2}
\begin{cases}
 { < 0,} & {x > 1,} \\
 { > 0,} & {x < 1.} \\
\end{cases}
\end{equation}

In view of (\ref{eq78}), we conclude that the function $g_{N_2 N_1
\_N_2 G} (x)$ is increasing in $x \in (0,1)$ and decreasing in $x
\in (1,\infty )$, and hence
\begin{equation} \label{eq79} M =
\mathop {\sup }\limits_{x \in (0,\infty )} g_{N_2 N_1 \_N_2 G} (x)
= g_{N_2 N_1 \_N_2 G} (1) = \frac{1}{3}.
\end{equation}

Applying the inequalities (\ref{eq34}) for the measures $M_{N_2
N_1 } (P\vert \vert Q)$ and $M_{N_2 G} (P\vert \vert Q)$ along
with (\ref{eq79}) we get the required result.
\end{proof}

\begin{proposition} \label{pro23} The following inequality hold:
\begin{equation}
\label{eq80}
M_{N_2 G} (P\vert \vert Q) \leqslant \frac{3}{4}M_{AG} (P\vert \vert Q).
\end{equation}
\end{proposition}

\begin{proof} Let us consider
\begin{equation}
\label{eq81} g_{N_2 G\_AG} (x) = \frac{{f}''_{N_2 G}
(x)}{{f}''_{AG} (x)} = - \frac{1 + x^{3 / 2} - (2x + 2)^{3 /
2}}{(2x + 2)^{3 / 2}}, \,\, x \in (0,\infty ),
\end{equation}

\noindent where ${f}''_{N_2 G} (x)$ and ${f}''_{AG} (x)$ are as
given by (\ref{eq68}) and (\ref{eq64}) respectively.

From (\ref{eq81}), we have
\begin{equation}
\label{eq82}
{g}'_{N_2 G\_AG} (x) = \frac{3\left( {1 - \sqrt x } \right)}{(2x + 2)^{5 /
2}}
\begin{cases}
 { \leqslant 0,} & {x \geqslant 1}, \\
 { \geqslant 0,} & {x \leqslant 1}. \\
\end{cases}
\end{equation}

In view of (\ref{eq82}), we conclude that the function $g_{AH\_N_2
N_1 } (x)$ is increasing in $x \in (0,1)$ and decreasing in $x \in
(1,\infty )$, and hence
\begin{equation}
\label{eq83}
M = \mathop {\sup }\limits_{x \in (0,\infty )} g_{N_2 G\_AG} (x) = g_{N_2
G\_AG} (1) = \frac{3}{4}.
\end{equation}

Applying the inequalities (\ref{eq34}) for the measures $M_{N_2 G}
(P\vert \vert Q)$ and $M_{AG} (P\vert \vert Q)$ along with
(\ref{eq83}) we get the required result.
\end{proof}

\begin{proposition} \label{pro24} The following inequality hold:
\begin{equation}
\label{eq84}
\frac{1}{4}M_{AG} (P\vert \vert Q) \leqslant M_{AN_2 } (P\vert \vert Q).
\end{equation}
\end{proposition}

\begin{proof} Let us consider
\begin{equation}
\label{eq85} g_{AG\_AN_2 } (x) = \frac{{f}''_{AG} (x)}{{f}''_{AN_2
} (x)} = \frac{(2x + 2)^{3 / 2}}{\left( {\sqrt x + 1}
\right)\left( {x - \sqrt x + 1} \right)}, \,\, x \in (0,\infty ),
\end{equation}

\noindent where ${f}''_{AG} (x)$ and ${f}''_{AN_2 } (x)$ are as
given by (\ref{eq64}) and (\ref{eq70}) respectively.

From (\ref{eq85}), we have
\begin{equation}
\label{eq86}
{g}'_{AG\_AN_2 } (x) = \frac{3\left( {1 - \sqrt x } \right)\sqrt {2x + 2}
}{\left( {\sqrt x + 1} \right)^2\left( {x - \sqrt x + 1} \right)^2}
\begin{cases}
 { \leqslant 0,} & {x \geqslant 1} \\
 { \geqslant 0,} & {x \leqslant 1} \\
\end{cases}.
\end{equation}

In view of (\ref{eq86}), we conclude that the function
$g_{AG\_AN_2 } (x)$ is increasing in $x \in (0,1)$ and decreasing
in $x \in (1,\infty )$, and hence
\begin{equation}
\label{eq87}
M = \mathop {\sup }\limits_{x \in (0,\infty )} g_{AG\_AN_2 } (x) =
g_{AG\_AN_2 } (1) = 4.
\end{equation}

Applying the inequalities (\ref{eq34}) for the measures $M_{AG}
(P\vert \vert Q)$ and $M_{AN_2 } (P\vert \vert Q)$ along with
(\ref{eq87}) we get the required result.
\end{proof}

Combining the results given in the Propositions
\ref{pro21}-\ref{pro24}, we get the proof of the theorem.

The expression (\ref{eq72}) can also be written as
\begin{equation}
\label{eq88} \frac{1}{16}\Delta (P\vert \vert Q) \leqslant M_{N_2
N_1 } (P\vert \vert Q) \leqslant \frac{1}{3}M_{N_2 G} (P\vert
\vert Q)
\end{equation}
\[
 \leqslant \frac{1}{4}h(P\vert \vert Q) \leqslant M_{AN_2 } (P\vert \vert
Q).
\]

\begin{remark}
\begin{itemize}
\item[(i)] The classical divergence measures $I(P\vert \vert Q)$
and $J(P\vert \vert Q)$ appearing in the Section 1 can be written
in terms of \textit{Kullback-Leibler's relative information} as
follows:
\begin{equation}
\label{eq90}
I(P\vert \vert Q) = \frac{1}{2}\left[ {K\left( {P\vert \vert \frac{P +
Q}{2}} \right) + K\left( {Q\vert \vert \frac{P + Q}{2}} \right)} \right]
\end{equation}

\noindent and
\begin{equation}
\label{eq91}
J(P\vert \vert Q) = K(P\vert \vert Q) + K(Q\vert \vert P).
\end{equation}

Also we can write
\begin{equation}
\label{eq92}
J(P\vert \vert Q) = 4\left[ {I(P\vert \vert Q) + T(P\vert \vert Q)}
\right],
\end{equation}

\noindent where
\begin{equation}
\label{eq93}
T(P\vert \vert Q) = \frac{1}{2}\left[ {K\left( {\frac{P + Q}{2}\vert \vert
P} \right) + K\left( {\frac{P + Q}{2}\vert \vert Q} \right)} \right]
\end{equation}
\[
 = \sum\limits_{i = 1}^n {A(p_i ,q_i )\ln \left( {\frac{A(p_i ,q_i )}{G(p_i
,q_i )}} \right)} ,
\]

\noindent is the \textit{arithmetic and geometric mean divergence
measure} due to Taneja \cite{tan3}.

\item[(ii)] Recently, Taneja \cite{tan8} proved the following
inequality:
\begin{equation}
\label{eq94} \frac{1}{4}\Delta (P\vert \vert Q) \leqslant I(P\vert
\vert Q) \leqslant h(P\vert \vert Q) \leqslant 4\,M_{AN_2 }
(P\vert \vert Q)
 \leqslant \frac{1}{8}J(P\vert \vert Q) \leqslant T(P\vert \vert Q).
\end{equation}

Following the lines of the Propositions \ref{pro21}-\ref{pro24},
we can also show that
\begin{equation}
\label{eq95} \frac{1}{4}I(P\vert \vert Q) \leqslant M_{N_2 N_1 }
(P\vert \vert Q).
\end{equation}

Thus combining (\ref{eq88}) with (\ref{eq92}), (\ref{eq94}) and
(\ref{eq95}), we get the following inequalities among the
\textit{classical }and \textit{mean divergence measures}:
\begin{align}
\label{eq96} & \frac{1}{16}\Delta (P\vert \vert Q) \leqslant
\frac{1}{4}I(P\vert \vert Q) \leqslant M_{N_2 N_1 } (P\vert \vert
Q)\\
 & \quad \leqslant \frac{1}{3}M_{N_2 G} (P\vert \vert Q) \leqslant
\frac{1}{4}h(P\vert \vert Q) \leqslant M_{AN_2 } (P\vert \vert
Q)\notag\\
&\quad \quad \leqslant \frac{1}{32}J(P\vert \vert Q) \leqslant
\frac{1}{4}T(P\vert \vert Q) \leqslant \frac{1}{16}J(P\vert \vert
Q).\notag
\end{align}
\end{itemize}
\end{remark}

\end{document}